\documentclass[12pt]{article}
\usepackage{leqno}
\usepackage{amssymb}
\usepackage{amsmath}

\begin{document}
\newcommand{\bfR}{\hbox{\bbbld R}}
\newcommand {\prob}{probability}
\newcommand{\vsp}{\vskip 1em}
\newcommand{\Vsp}{\vskip 2em}
\newcommand{\bea}{\begin{eqnarray}}
\newcommand{\eea}{\end{eqnarray}}
\def\Sym{\mathop{\rm Sym}}
\newcommand{\suchthat}{\mid}
\newcommand{\halo}[1]{\Int(#1)}
\def\Int{\mathop{\rm Int}}
\def\Re{\mathop{\rm Re}}
\def\Im{\mathop{\rm Im}}
\def \qed {\hfill \vrule height6pt width 6pt depth 0pt}
\newcommand{\union}{\cup}
\newcommand{\goesto}{\rightarrow}
\newcommand{\bdy}{\partial}
\newcommand{\n}{\noindent}
\newcommand{\ds}{\displaystyle}
\newcommand{\be}{\begin{equation}}
\newcommand{\ee}{\end{equation}}
\newcommand{\ben}{\begin{eqnarray*}}
\newcommand{\een}{\end{eqnarray*}}
\newtheorem{theorem}{Theorem}[section]
\newtheorem{assertion}{Assertion}[section]
\newtheorem{proposition}{Proposition}[section]
\newtheorem{remark}{Remark}[section]
\newtheorem{lemma}{Lemma}[section]
\newtheorem{definition}{Definition}[section]
\newtheorem{claim}{Claim}[section]
\newtheorem{corollary}{Corollary}[section]
\newtheorem{observation}{Observation}[section]
\newtheorem{conjecture}{Conjecture}[section]
\newtheorem{question}{Question}[section]
\newtheorem{example}{Example}[section]
\newbox\qedbox
\newenvironment{proof}{\smallskip\noindent{\bf Proof.}\hskip \labelsep}%
                                {\hfill\penalty10000\copy\qedbox\par\medskip}
\newenvironment{proofspec}[1]%
                      {\smallskip\noindent{\bf Proof of #1.}\hskip \labelsep}%
 {\nobreak\hfill\hfill\nobreak\copy\qedbox\par\medskip}
\newenvironment{acknowledgements}{\smallskip\noindent{\bf Acknowledgements.}%
        \hskip\labelsep}{}

\newcounter{num}
\newcounter{alf}
\newcounter{rom}
\newcounter{Alf}

\title{The Elementary Polynomials in Noncommuting Variables}
\author{Samuel S. Holland, Jr.\\
Department of Mathematics and Statistics,\\ 
University of Massachusetts\\
 Box 34515\\
Amherst, MA 01003-4515\\
email:holland@@math.umass.edu}
\date{}
\thispagestyle{empty}

\maketitle

\begin{abstract}
We study the ring generated over a field of characteristic $0$ by
noncommuting indeterminates $\{x_{1},x_{2},...,x_{n}\}$ subject only
to the relations $x_{i}\sigma_{k}=\sigma_{k}x_{i}\;,\;i,k=1,2,...,n,$
and their consequences, where
$\sigma_{k}=\sigma_{k}(x_{1},x_{2},...,x_{n})$ is the k-th elementary
polynomial in the noncommuting variables $x_{i}$.  We assume $n\geq 3$ throughout.
\end{abstract}

\n{\bf 1. BASICS}\\

\n Let $K$ represent a commutative field of characteristic $0$, let
$x_{1},x_{2},...,x_{n}$ be algebraically independent noncommuting
indeterminates over $K$, and denote by $P=K<x_{1},x_{2},...,x_{n}>$ the
free $K$-ring generated over $K$ by the $x_{i}$.  Each $x_{i}$
commutes with every element of $K$. {\em Monomials}  $u\in P$ have the
form

\begin{equation}
u=x^{j(1)}_{i(1)}x^{j(2)}_{i(2)}\cdots x^{j(\ell)}_{i(\ell)},
\end{equation}

\n and will be understood to satisfy the following conditions: the
$j(\cdot)$ are integers $\geq 1, \; i(\cdot)\in\{1,2,...,n\},$ and
$i(p)\neq i(q)$ when $|p-q|=1$. We call $\{j(1),j(2),...,j(\ell)\}$
the {\em exponent sequence} of $u,$ the sequence of subscripts
$\{i(1),i(2),...,i(\ell)\}$ its {\em complexion}, $j(1)+j(2)+\cdots
+j(\ell)$ its {\em degree}, and $\ell$ its {\em length}.  We denote
the set of all monomials by $Q; Q\cup\{1\}$ is a vector space basis of
$P$ over $K$.

We define the {\em elementary polynomials,}
$\sigma_{k}(x_{1},x_{2},...,x_{n}),$ as follows.  Set 
$\sigma_{0}=1,\;\;\sigma_{k}=0$ when $k<0\;\;$ or $\;\;k>n,$ and, for,
$k=1,2,...,n,$

\begin{equation}
\begin{array}{l}
\sigma_{k}(x_{1},x_{2},...,x_{n})=\sum(x_{i(1)}x_{i(2)}\cdots
x_{i(k)}:\;i(\cdot)\in\{1,2,...,n\},\\
{\makebox{\hspace*{2in}}} i(1)<i(2)<...<i(k))
\end{array}
\end{equation}

\n Each $\sigma_{k},\;\;k=1,2,...,n$, has $\begin{pmatrix} 
n\\
k
\end{pmatrix}$ monomials where $\begin{pmatrix}
n\\
k
\end{pmatrix}=n!/(n-k)!k!$ is the usual binomial coefficient. We may
break up the sum (2) into $n-k+1$ subsums according as the
monomials involved begin with $x_{1},x_{2},...$.

\begin{equation}
\left\{\begin{array}{ll}
\sigma_{k}=\sum (x_{1}x_{i(2)}\cdots x_{i(k)}:1<i(2)<...<i(k)) &
\left(\begin{array}{l} n-1\\
k-1
\end{array}\right) \text{nonomials}\\
+\sum (x_{2}x_{i(2)}\cdots x_{i(k)}:2 < i(2)...< i (k)) &
\left(\begin{array}{l}
n-2\\
k-1
\end{array}\right) \text{monomials}\\
+\cdots & \\
+ x_{n-k+1} x_{n-k+2}\cdots x_{n} & \left(\begin{array}{l}
k-1\\
k-1
\end{array}\right)= 1\; \text{monomial}
\end{array}\right.
\end{equation}

\n In more traditional form,

\[\begin{array}{l}
 \sigma_{1}(x_{1},x_{2},...,x_{n})=x_{1}+x_{2}+\cdots+x_{n}\\
\sigma_{2}(x_{1},x_{2}...,x_{n})=(x_{1}x_{2}+x_{1}x_{3}+\cdots+x_{1}x_{n})+(x_{2}x_{3}+\cdots
+x_{2}x_{n})+\cdots + x_{n-1}x_{n}\\
{\makebox{\hspace*{2.5in}}}......\\
\sigma_{n-1}(x_{1},x_{2},...,x_{n})=(x_{1}x_{2}\cdots
x_{n-1}+x_{1}x_{2}\cdots x_{n-2}x_{n}+\cdots
+x_{1}x_{3}\cdots x_{n})+x_{2}x_{3}\cdots x_{n}\\
\sigma_{n}(x_{1},x_{2},... ,x_{n})=x_{1}x_{2}\cdots x_{n}
\end{array}\]

\n In the commutative case, these polynomials are the same as
the ``elementary symmetric functions'' (as commonly written), but here
the variables do not commute so specification of the order of the
variables is necessary.

We may rewrite equation (3) as follows,

\[\begin{array}{lcl}
\sigma_{k}(x_{1},x_{2},...,x_{n}) & = & \sum\;(x_{1}x_{i(2)}\cdots
x_{i(k)}: \;i(\cdot)\in\{2,3,...,n\},\\
 & & {\makebox{\hspace*{1.2in}}} i(2)<i(3)<...<i(k))\\
& + & \sum\; (x_{i(1)}x_{i(2)}\cdots
x_{i(k)}:\;i(\cdot)\in\{2,3,...,n\},\\
& & {\makebox{\hspace*{1.2in}}}i(1)<i(2)<...<i(k))
\end{array}\]

\n and deduce therefrom the {\em first recursion formula} for
$\sigma_{k}$\\

\n(I)\;\;\;$\displaystyle{\begin{array}{l}
\sigma_{k}(x_{1},x_{2},...,x_{n},)=x_{1}\sigma_{k-1}(x_{2},x_{3},...,x_{n})+\sigma_{k}(x_{2},x_{3},...,x_{n}),\\
 {\makebox{\hspace*{3.5in}}}k=1,2,...n
\end{array}}$\\

\n which is valid over the $k$-range indicated because we have
set $\sigma_{0}=1$ and $\sigma_{k}=0$ when $k>\#$ of variables.

We may derive a companion formula to (2) where we factor out the last
variable $x_{n}$ instead of the first, and, doing so, we get the {\em
  second recursion formula} for $\sigma_{k}$.\\

\n(II)\;\;\;$\displaystyle{\begin{array}{l}
\sigma_{k}(x_{1},x_{2},...,x_{n})=\sigma_{k-1}(x_{1},x_{2},...,x_{n-1})x_{n}+\sigma_{k}(x_{1},x_{2},...,x_{n-1}),\\
{\makebox{\hspace*{3.5in}}}k=1,2,...,n.
\end{array}}$\\

\n We are concerned here with the action on the free ring $P$
of the group $G$ of circular permutations by which I mean the cyclic
subgroup of the full permutation group $S_{n}$ on $n$ letters
generated by $(1 2...n)$.  The action of $g\in G$ on the monomial (1)
is given by

\[u^{g}=x^{j(1)}_{i(i)g}x^{j(2)}_{i(2)g}\cdots
x^{j(\ell)}_{i(\ell)g}.\]

\n Note that $G$ acts only on the complexion of $u$; its
exponent sequence does not change. The action of $G$ extends by
linearity to polynomials.  The coefficient field $K$ is pointwise fixed
by $G$.  A polynomial $p\in P$ is {\em invariant under $G$} when
$p^{g}=p$ for all $g\in G$.  As $G$ is cyclic, it is enough that
$p^{g(1)}=p$ for the generator $g(1)=(1 2...n)$ of $G$.\\

\n{\bf Theorem 1.1.} {\em In the free ring $P$, the (two-sided)
  ideal generated by the
  $\sigma_{k}(x_{1},x_{2},...,x_{n})-\sigma_{k}(x_{1},x_{2},...,x_{n})^{g}\;,\; k=1,2,...,n,\forall g\in G$, is the same as that generated by the commutators $[x_{i},\sigma_{k}]=x_{i}\sigma_{k}(x_{1},x_{2},...,x_{n})-\sigma_{k}(x_{1},x_{2},...,\\
x_{n})x_{i},1\leq \imath,k\leq n$.  Or, what amounts to the same thing, the set of relations $\{\sigma^{g}_{k}=\sigma_{k}:1\leq k\leq n,\;\forall g\in G\}$ is equivalent to the set $\{x_{i}\sigma_{k}=\sigma_{k}x_{i}:\;1\leq i, k\leq n\}.$}\\

\n{\bf Proof.} The following notations will be useful.  Denote
by $g(1)=(1 2 ...n)$ the generator of the group $G$ of circular
permutations, and by $g(i)$ the $i$-th power of $g(1)$.  $S\circ$
$\;\;G=\{1=g(0), g(1),g(2),..., g(n-1)\}$. The recursion formulas (I)
and (II) hold for any set of $n$ subscripted letters as long as their
order there is maintained.  In particular, we may apply $g(i)$ to them
to get \\

\n
I(i)$\displaystyle{\;\;\;\;\;\;\;\;\;\;\sigma_{k}^{g(i)}=x_{i+1}\sigma_{k-1}(x_{i+2},...,x_{i})+\sigma_{k}(x_{i+2},...,x_{i})}$\\

\n
I(i+1)$\displaystyle{\;\;\;\;\;\sigma_{k}^{g(i+1)}=\sigma_{k-1}(x_{i+2},...,x_{i})x_{i+1}+\sigma_{k}(x_{i+2},...,x_{i})}$\\

\n for $k=1,2,...,n$ and $i=0,1,...,n-1$.  (In the complexions
we work mod n.)  Denote by $I$ the ideal generated by the
$\sigma_{k}-\sigma_{k}^{g}$, and by $I^{\prime}$ that generated by
the $[x_{i},\sigma_{k}]$. We aim to show that $I=I^{\prime}$. (In this
paper ``ideal'' means two-sided ideal.) \\

\n $\displaystyle{I\subset}$ $I^{\prime}:$ The recursion I$(0)$
  gives

\[\sigma_{k}=x_{1}\sigma_{k-1}(x_{2},...,x_{n})+\sigma_{k}(x_{2},...,x_{n}),k=1,2,...,n\]

\n $\displaystyle{\begin{array}{lcl}\text{so}\;\;\;\;[x_{1},\sigma_{k}] &=
  &[x_{1},x_{1}\sigma_{k-1}(x_{2},...,x_{n})]+[x_{1},\sigma_{k}(x_{2},...,x_{n})]\\
& = &
x_{1}[x_{1},\sigma_{k-1}(x_{2},...,x_{n})]+[x_{1},\sigma_{k}(x_{2},...,x_{n})]\end{array}}$\\

\n (using $[a,bc]=[a,b]c+b[a,c])$. Put $k=1$. As $\sigma_{0}=1$
we get

\[[x_{1},\sigma_{1}]=[x_{1},\sigma_{1}(x_{2},...,x_{n})]\]

\n so the latter commutator belongs to $I^{\prime}$.  Put
$k=2$.  

\[[x_{1},\sigma_{2}]=x_{1}[x_{1},\sigma_{1}(x_{2},...,x_{n})]+[x_{1},\sigma_{2}(x_{2},...,x_{n})].\]

\n We have just shown that the first term on the right belongs
to $I^{\prime}$, hence also $[x_{1},\sigma_{2}(x_{2},...,x_{n})]$.
Continuing, we establish that

\[[x_{1},\sigma_{k}(x_{2},...,x_{n})]\in
I^{\prime},\;\;\;k=1,2,...,n\]

\n (The case $k=n$ can be included because $\sigma_{k}=0$ when
$k>$\# of variables.) Now go to II$(1)$,

\[\sigma^{g(1)}_{k}=\sigma_{k-1}(x_{2},...,x_{n})x_{1}+\sigma_{k}(x_{2},...,x_{n}),\]

\n and subtract I$(0)$-II$(1)$,

\[\sigma_{k}-\sigma_{k}^{g(1)}=[x_{1},\sigma_{k-1}(x_{2},...,x_{n})]\in
I^{\prime},\;\;k=1,2,...,n.\]

\n That is the first step.  If we had taken the trouble to
establish in advance that $(I^{\prime})^{g}\subset I^{\prime},$ we
  would be done.  But this direct method seems preferable, so, we
  proceed to the second step.  From $\sigma_{k}-\sigma_{k}^{g(1)}\in
  I^{\prime}$ we get
  $[x_{i},\sigma_{k}-\sigma_{k}^{g(1)}]=[x_{i},\sigma_{k}]-[x_{i},\sigma_{k}^{g(1)}]\in I^{\prime}$, hence $[x_{i},\sigma_{k}^{g(1)}]\in I^{\prime},\;\;i,k=1,2,...,n$.  The recursion I$(1)$ gives

\[\sigma_{k}^{g(1)}=x_{2}\sigma_{k-1}(x_{3},...,x_{1})+\sigma_{k}(x_{3},...,x_{1})\;\;\;\text{so}\]

\[[x_{2},\sigma_{k}^{g(1)}]=x_{2}[x_{2},\sigma_{k-1}(x_{3},...,x_{1}]+[x_{2},\sigma_{k}(x_{3},...,x_{1})].\]

\n Put $k=1,k=2$, etc, as before, to show that
$[x_{2},\sigma_{k}(x_{3},...,x_{1})]\in
I^{\prime},\;k=1,2,...,n$. Then subtract I$(1)-$II$(2),$

\[\sigma_{k}^{g(1)}-\sigma^{g(2)}_{k}=[x_{2},\sigma_{k-1}(x_{3},...,x_{1})]\in
I^{\prime},\;\;k=1,2,...,n,\]

\n which completes the second step.  Continue, to get

\[\sigma_{k}^{g(i)}-\sigma_{k}^{g(i+1)}\in
I^{\prime},\;\;k=1,2,..,n,\;\;i=0,1,...,n-1\]

\n which shows $I\subseteq I^{\prime}$\\

\n $I^{\prime}\subseteq I$. We need to show that all commutators
$[x_{i},\sigma_{k}]$ belong to the ideal I generated by the
$\sigma_{k}-\sigma_{k}^{g}$. Put $i=0$ in I(i) and II(i+1), \\

\n
I(0)$\displaystyle{\;\;\;\;\;\;\;\;\sigma_{k}=x_{1}\sigma_{k-1}(x_{2},...,x_{n})+\sigma_{k}(x_{2},...,x_{n})}$\\

\n
II(1)\;\;\;\;\;\;\;$\displaystyle{\sigma_{k}^{g(1)}=\sigma_{k-1}(x_{2},..,x_{n})x_{1}+\sigma_{k}(x_{2},...,x_{n})}$\\

\n and subtract to get

\[\sigma_{k}-\sigma_{k}^{g(1)}=[x_{1},\sigma_{k-1}(x_{2},...,x_{n})],\;\;k=1,2,...,n.\]

\n Then, using I(0),

\[[x_{1},\sigma_{k}]=x_{1}(\sigma_{k}-\sigma_{k}^{g(1)})+(\sigma_{k+1}-\sigma_{k+1}^{g(1)}),k=1,2,...,n,\]

\n which shows that $[x_{1},\sigma_{k}]\in I, k=1,2,...,n.$
Next, put $i=1$ in I(i) and II(i+1)\\

\n
I(1)\;\;\;\;\;\;\;\;\;\;$\displaystyle{\sigma_{k}^{g(1)}=x_{2}\sigma_{k-1}(x_{3},...,x_{1})+\sigma_{k}(x_{3},...,x_{1})}$\\

\n
II(2)\;\;\;\;\;\;\;\;\;\;$\displaystyle{\sigma_{k}^{g(2)}=\sigma_{k-1}(x_{3},...,x_{1})x_{2}+\sigma_{k}(x_{3},...,x_{1}),}$\\

\n and subtract to get

\[\sigma_{k}^{g(1)}-\sigma_{k}^{g(2)}=[x_{2},\sigma_{k-1}(x_{3},...,x_{1})]\;\;k=1,2,...,n.\]

\n Then, using I(1),

\[[x_{2},\sigma_{k}^{g(1)}]=x_{2}(\sigma_{k}^{g(1)}-\sigma_{k}^{g(2)})+(\sigma_{k+1}^{g(1)}-\sigma_{k+1}^{g(2)}).\]

\n Then write
$\sigma_{k}=\sigma_{k}^{g(1)}+(\sigma_{k}-\sigma_{k}^{g(1)})$ to get
$[x_{2},\sigma_{k}]=x_{2}(\sigma_{k}^{g(1)}-\sigma_{k}^{g(2)})+(\sigma_{k+1}^{g(1)}-\sigma_{k+1}^{g(2)})+[x_{2},(\sigma_{k}-\sigma_{k}^{g(1)})]\;\;k=1,2,...,n.$
which shows that $[x_{2},\sigma_{k}]\in I,\;\;k=1,2,...n.$\\

\n The remaining steps follow the same pattern.  That completes
the proof of Theorem 1.1.\\

\n The proof of Theorem 1.1 applies without change to the free
$\Bbb{Z}$-ring ${\Bbb{Z}}<x_{1},x_{2},...,x_{n}>$ where ${\Bbb{Z}}\subseteq$
$K$ is the ring of rational integers.  As no constant term is involved,
we have\\

\n{\bf Corollary 1.2} {\em Let $B$ be any ring,
  $b_{1},b_{2},...,b_{n}$ any finite number of elements from $B$.  The
  following two conditions are equivalent.}

\begin{list}%
{(\roman{rom})}{\usecounter{rom}\setlength{\rightmargin}{\leftmargin}}
{\em\item Each $\sigma_{k}(b_{1},b_{2},...,b_{n}),k=1,2,...,n,$ is
  invariant under circular permutations of the $b_{i}$.
\item Each $\sigma_{k}(b_{1},b_{2},,,b_{n})$ commutes with each
$b_{i},i,k=1,2,...,n.$}
\end{list}

\n We shall denote by $I$ the common ideal described in Theorem
1.1, and shall denote by $R^{\prime}$ the quotient ring $P/I$.  We
shall work with $R^{\prime}$ as a polynomial ring with relations
(those specified in Theorem 1.1). We shall also make frequent use of
the fact that a polynomial in $R^{\prime}$ is zero there if and only
if it, exactly as written, when considered in the free ring $P$,
belongs to the ideal $I$.  Because a polynomial in $P$ belongs to $I$
if and only if its homogeneous components do, a polynomial in
$R^{\prime}$ will be zero there if and only if its homogeneous
components are separately zero. Thus two polynomials in $R^{\prime}$
will be equal exactly when their homogenous components are separately equal, and monomials of different degree
cannot be equal.  But monomials of the same degree can be; for example
when $n=3$
$x_{1}x_{2},x_{3}=x_{2}x_{3}x_{1}=x_{3}x_{1}x_{2}(=\sigma_{3})$.
Further the same element of $R^{\prime}$ may have different polynomial
representations.  Again, when $n=3, x_{1}x_{2}+x_{1}x_{3}+x_{2}x_{3}=x_{2}x_{3}+x_{2}x_{1}+x_{3}x_{1}(=\sigma_{2})$ so,
in $R^{\prime}, x_{1}x_{2}+x_{1}x_{3}-x_{2}x_{1}-x_{3}x_{1}=0$.\\

\n The general homogeneous polynomial of degree $w$ in $I$ is a
$K$-linear combination of polynomials $u[x_{i},\sigma_{k}]v$ where $u$
and $v$ are monomials of degree $r$ and $s$ respectively (say) and
$w=r+s+k+1$.  As $k\geq 1,w\geq 2$, so $I$ contains no linear
polynomials.  Hence, in $R^{\prime}$, the variables $x_{i}$ are
$K$-linearly independent.  If $w=2$, then $r=s=0$ and $k=1$ in which
case the generators $[x_{i},\sigma_{k}]$ are just sums of the
$[x_{i},x_{j}].$ We shall use the notations
$[i,j]=x_{i}x_{j}-x_{j}x_{i}$ for these additive commutators, and we
need only consider $[i,j]$ when $i<j$. We shall refer to commutators
$[i,i+1]$ as ``diagonal'', and $[i,j]$ when $j>i+1$ as
``off-diagonal''.
\newpage

\n{\bf Theorem 1.3} {\em The off-diagonal commutators $[i,j+1],
  1\leq i < j< n$ are $K$-linearly independent in $R^{\prime}$.  The
  diagonal commutators are expressed in terms of the off-diagonal by
  the formulas

\begin{equation}
[k,k-1]=\sum^{k-2}_{p=1}\sum^{n}_{j=k}\;[p,j]-\sum^{n}_{j=k+1}\;[k-1,j],
\end{equation}

\n $k=2,3,...,n,$, where we use the convention that sums whose
lower limit exceeds the upper are given the value $0$.}\\

\n Under the stated convention, the first sum vanishes when
$k=2$, and the last does when $k=n$.  Hence, when $k=2$, formula (4)
becomes

\[[2,1]=\sum^{n}_{j=3}\;[1,j],\]

\n and, when $k=n$,

\[[n,n-1]=\sum^{n-2}_{p=1}\;[p,n].\]

\n Note also this consequence when $n=3:
[1,2]=[2,3]=[3,1]$. (Remember $n\geq 3$ is assumed throughout.)\\

\n{\bf Proof:} As commutators are homogeneous of degree 2, the
assertion that a $K$-linear combination of off-diagonal commutators
equals zero in $R^{\prime}$ takes the form

\[\sum_{1\leq i< j<
  n}\;a_{ij}[i,j+1]=\sum^{n}_{i=1}\;b_{i}(x_{i}\sigma_{1}-\sigma_{1}x_{i})\]

\n in the free ring $P$, where the a's and b's belong to $K$.
An easy computation evaluates the right side and we get

\begin{equation}
\sum_{1\leq
  i<j<n}\;a_{ij}[i,j+1]=\sum^{n}_{i=1}\;b_{i}\;\sum^{n}_{\ell=i}\;[i,\ell].
\end{equation}

\n Because the left side of (5) consists solely of off-diagonal
commutators, it does not contain either $[1,2]$ or $[2,1]$.  The right
side however contains $b_{2}[2,1]+b_{1}[1,2].$ As we are working in a
free ring we must have $b_{1}=b_{2}$. Arguing similarly with $[2,3]$
we get $b_{2}=b_{3}$.  And so on.  The right side of (5) thus becomes

\[b\sum^{n}_{i=1}\;\sum^{n}_{\ell=1}\;[i,\ell]=0\]

\n which is zero because we are summing over all entries of a
skew symmetric matrix.  Thus the right side of (5) is zero in the free
ring $P$.  And, as we are working in the free ring, it follows that all
$a_{ij}$ must be zero.  This proves the linear independence of the
off-diagonal commutators in $R^{\prime}$.  Next we turn to equation
(4) expressing the diagonal commutators in terms of the
off-diagonal. In the ring $R^{\prime}$ we have the identity

\[\sum^{n}_{p=1}\;a_{p}(x_{p}\sigma_{1}-\sigma_{1}x_{p})=0\]

\n for any selection of $a_{p}\in K$, because $x\sigma=\sigma x$
holds identically in $R^{\prime}$.  The above equation may be recast
in the form

\begin{equation}
\sum^{n}_{1\leq p< q\leq n}\;(a_{p}-a_{q})[p,q]=0.
\end{equation}

\n In (6) , put $a_{1}=1$, and $a_{p}=0,p\geq 2$ to get
$[1,2]+[1,3]+\cdots+[1,n]=0$ which we may solve for $[1,2]$ in terms
of off-diagonal commutators.  Next put
$a_{1}=a_{2}=1,\;\;a_{p}=0,p\geq 3$ and solve for $[2,3]$.  Continue,
to get $n-1$ equations for the $[k-1,k]$ which constitute the system
of equations (4). We have $n-1$ independent equations for $n(n-1)/2$
unknowns, and we have solved for the $n-1$ diagonal commutators in
terms of the $(n-1)(n-2)/2$ independent off-diagonal ones.  That
completes the proof of Theorem 1.3.\\

\n {\bf Corollary 1.4.} {\em In $R^{\prime},\;\;x_{i}x_{j}\neq
  x_{j}x_{i}$ when $i\neq j$.}\\

\n {\bf Proof.} As the off-diagonal commutators are
$K$-linearly independent, none of them can be zero.  And formula (4)
shows that no diagonal commutators can be zero either lest we get a
nontrivial linear combination of the off-diagonals equal to
zero. $\diamond$\\

\n We make no use in the sequel of the following corollary, so
I register it here without proof (the ``$j$'' in front of the bracket is
the integer $j$).\\

\n {\bf Corollary 1.5.} {\em For the sum of the diagonal
  commutators in $R^{\prime}$ we have}

\[\sum^{n}_{k=2}\;[k,k-1]=\sum^{n-2}_{p=1}\sum^{n-p}_{j=2}\;j[p,j+p].\]

\n The next corollary that describes a canonical form for the
homogeneous quadratic in $R^{\prime}$ can be deduced from Theorem 1.3
in a straightforward way, so that proof is left to the reader.\\

\n {\bf Corollary 1.6.} {\em Every homogeneous quadratic
  polynomial $p\in R^{\prime}$ may be written

\[p=\sum^{n}_{i=1}\;a_{i}x_{i}^{2}+\sum_{1\leq i<j\leq
  n}\;b_{ij}x_{i}x_{j}+\sum_{1\leq i<j < n}\;c_{ij}[i,j+1]\]

\n where the coefficients (in $K$) are uniquely determined by
$p$.}\\

\n Our construction of the ring $R^{\prime}$ can be viewed as
parallel to a construction of the free commutative ring.  One may show
that, in the free ring $P=K<x_{1},x_{2},...,x_{n}>,$ the set of
relations $\{x_{i}x_{j}=x_{j}x_{i}:\;i,j=1,2,...,n\}$ is equivalent to
the set $\{\sigma_{k}=\sigma_{k}^{\pi}:\;\forall\pi\in S_{n}\}$.
Our Theorem 1.1 shows that if we restrict to invariance under only the
subgroup $G=<(1 2...n)>$ of $S_{n}$, then the corresponding
commutativity is $x_{i}\sigma_{k}=\sigma_{k}x_{i}$.  Each commutator
$[x_{i},\sigma_{k}]$ can be expressed in terms of the $[x_{i},x_{j}]$,
so the ideal $J$ of the free ring $P$ generated by the $[x_{i},x_{j}]$
contains $I$, and contains it properly when $n\geq 3$ by Corollary 1.4. (When
  $n=2$, $J=I$ which is why we assume $n\geq 3$ thoughout.) Hence, if
  $A=K[x_{1},x_{2},...,x_{n}]$ is the free abelian ring in $n$
  commuting indeterminates, we have

\[ A\cong P/J\cong(P/I)/(J/I)\]

\n or

\[A\cong R^{\prime}/(J/I)\]

\n and thus the following diagram where $\omega$ denotes the
homomorphism $P\rightarrow P/J\cong A,\;\;\varphi:P\rightarrow  P/I\cong
R^{\prime},$ and $\Psi: R^{\prime}\rightarrow R^{\prime}/(J/I)\cong
A$.

\[
  \begin{picture}(270,130)(0,15)
    \put(55,135){$P =K<x_{1}x_{2},...,x_{n}>$ the free ring}
     \put(70,130){\vector(0,-1){95}}
     \put(50,130){\vector(-1,-1){35}}
    \put( 75, 80){$\omega=\Psi \circ\varphi$}
    \put(55, 20){$A =K[x_{1},x_{2},...,x_{n}]$ the free abelian ring}
    \put(  0, 80){$R^{\prime}$}
     \put(15, 75){\vector(1,-1){40}}
    \put( 27,120){$\varphi$}
    \put( 27, 43){$\Psi$}
  \end{picture}
\]
\begin{center}
{Figure 1}
\end{center}

\noindent{\bf Theorem 1.7} 
{\em \begin{list}%
{(\roman{rom})}{\usecounter{rom}\setlength{\rightmargin}{\leftmargin}}
\item A necessary condition that $p(x_{1},x_{2},...,x_{n})=0$ in
  $R^{\prime}$ is that it can be brought to zero by assuming that the
  $x_{i}$ mutually commute.
\item In $R^{\prime}$ the elementary polynomials
  $\sigma_{k}(x_{1},x_{2},...,x_{n}),k=1,2,..,n$ are algebraically
  independent over $K$, thus generate a commutative domain
  $W^{\prime}=K[\sigma_{1},\sigma_{2},...,\sigma_{n}]$ within the
  center of $R^{\prime}$.
\end{list}}

\n{\bf Proof.} As for (i), if $p=0$ in $R^{\prime}$, then
$\Psi(p)=0$ in the free abelian ring $A$.  So if $p$ cannot be brought
to zero by assuming that its variables commute, then it cannot be zero
in $R^{\prime}$. As for (ii), note first that the $\sigma_{k}$ are
central in $R^{\prime}$ by construction.  If
$p(\sigma_{1},\sigma_{2},...,\sigma_{n})\in R^{\prime}$ is a
$K$-coefficient polynomial in the $\sigma_{k}$ (nontrivial), then
$\Psi(p)$ is the same polynomial which now lies in the free abelian
ring $A=K[x_{1},x_{2},...,x_{n}]$ where we know the ``elementary
symmetric functions'' $\sigma_{k}$ are algebraically independent.  So
$\Psi (p)\neq 0$ in $A$, hence $p\neq 0$ in $R^{\prime}$. $\diamond$\\

\n Look for a moment at the commutative case.  The commutative
domain $A=K[x_{1},x_{2},...,x_{n}]$ is contained in the field
$F=K(x_{1},x_{2},...,x_{n})$ of rational functions in the $x_{i}$.
The field $F$ contains the field
$W=K(\sigma_{1},\sigma_{2},...,\sigma_{n})$ of rational functions in
the algebraically independent elementary symmetric functions
$\sigma_{k}$, and $F$ is a Galois extension of $W$ with Galois group
$S_{n}$. $F/W$ is the splitting field of
$y^{n}-\sigma_{1}y^{n-1}+...+(-1)^{n}\sigma_{n}$. What follows is a
noncommutative analogue of this set-up.\\

\n{\bf Definition 1.8.} {\em By the symbol $R$ we shall denote
  the extension of $R^{\prime}$ constructed as described below.  $R$
  contains $R^{\prime}$ as a subring, and contains the field
  $W=K(\sigma_{1},\sigma_{2},...,\sigma_{n})$ within its center.}\\

\n The method of construction of the field of fractions of a
commutative domain works without change to produce the ring $R$
described above.  $R$ consists of equivalence classes of pairs
$(r,\lambda), r\in R^{\prime}, 0\neq\lambda\in
  W=K[\sigma_{1},\sigma_{2},...,\sigma_{n}]$ under the equivalence
  relation $(r_{1},\lambda_{1})=(r_{2},\lambda_{2})\Longleftrightarrow
  r_{1}\lambda_{2}=r_{2}\lambda_{1}$. The ring operations are defined
  as usual, and the map $r\mapsto (r,1)$ embeds $R^{\prime}$ in $R$.
  In working with $R$ we shall write $\sigma^{-1}$ for $(1,\sigma)$,
  etc.\\

\n The ring $R$ is the focus of this paper.  We note the
following easily verified facts about $R$: (1) Every $x_{i}\in R$ is a
root of the polynomial
$f(y)=y^{n}-\sigma_{1}y^{n-1}+\sigma_{2}y^{n-2}-\cdots
+(-1)^{n}\sigma_{n}$ ($y$ a central indeterminate) which factors in
$R$ as $f(y)=(y-x_{1})(y-x_{2})\cdots (y-x_{n})$ where the $x_{i}$ may
be circularly permuted.  (2) Every $x_{i}$ has an inverse in $R$ given
by
$x^{-1}=(-1)^{n+1}\sigma_{n}^{-1}(x^{n-1}-\sigma_{1}x^{n-2}+\cdots+(-1)^{n-1}\sigma_{n-1})(x=x_{i})$,
and $x_{i}^{-1}$ is also a root of an nth degree polynomial, coefficients
in $W$. (3) Every monomial has an inverse in $R$.\\

\n I cannot however settle this basic question: Is $R$ a
domain? I have very little to offer on either side of this question,
but, hoping for an affirmative answer, I state it as\\

\n {\bf Conjecture 1.9.} {\em $R$ is a domain.}\\

\n A key would seem to be to establish whether or not a product
of commutators $[x_{i},x_{j}], i\neq j,$ can be zero in $R$.\\

\vspace*{.2in}
\n {\bf 2. INVARIANCE UNDER  $G=<(1 2... n)>$}\\

\n In constructing the ring $R$ we have forced the elementary
polynomials $\sigma_{k}$ to be fixed under the group $G$ of circular
permutations.  $R$ will contain other polynomials fixed by $G$, namely
the images in $R^{\prime}=P/I$ of those polynomials in the free ring
$P$ already fixed by $G$ (because $I^{g}=I$ for all $g$ in $G$).  So
any information about the fixed ring of $P$ should, via the map
$\varphi: P\longrightarrow R^{\prime}\subseteq R$ provide information
about $R$. In [1] Bergman and
Cohn have solved a much more general problem.  They consider the free
product of a family of copies of a ring over a skew field, and
determine the structure of the subring fixed under the action of any
group that acts with finite orbits.  I am very grateful to
Prof. Bergman for bringing this paper to my attention, and for
providing valuable advice for its application to this particular case.
In this section I will effect the specialization of [1] to our
situation, the determination of the structure of the subring of the
free product of $n$ copies of $K[x]$ fixed under the action of the
group of circular permutations. For proofs, refer to
[1] and the paper [2] by Wolf.\\

\n We denote by $S$ the subring of the free ring $P$ consisting
of polynomials fixed by $G$.  Given a monomial $u\in P$, the {\em
  orbit polynomial} $\overline{u}$ of $u$ is

\[\overline{u}=\sum (u^{g}:g\in G),\]

\n and consists of exactly $n$ different monomials.  The orbit
polynomials, together with 1, consitute a vector-space basis of $S$
over $K$. If $\overline{u}$ and $\overline{v}$ are two orbit
polynomials, their product $\overline{u}$ $\overline{v}$ will also be
invariant under $G$ so will be a unique $K$-linear combination of
orbit polynomials.  That $K$-linear combination can be explicitly
described: In the term-by-term product $\overline{u}$ $\overline{v}$, all
$n^{2}$ monomials $u^{g}v^{h}$ are different, and when rearranged so
that each grouping consists of monomials carried into each other by
the action of $G$, each grouping is an orbit polynomial, and their sum
(all coefficients 1) is the unique $K$-linear representation of
$\overline{u}$ $\overline{v}$.\\

\n Next we introduce an ordering that is the key to proving that the
fixed ring $S$ is freely generated by certain explicitly given orbit
polynomials.  Bergman and Cohn, dealing with their more general
situation described earlier, order the complexions (finite
nonrepetitive sequences from the index set, here $\{1,2,...,n\}$).
For the polynomial ring under consideration here, I found it more
convenient to order the monomials directly, rather than their
complexions, using the first ordering given by Wolf [2; Sec. 3].  Her
ordering has all the essential properties of Bergman and Cohn's
``orbital odering'' [1; pp 526-527], and works perfectly well for the
group of circular permutations even though applied originally to the
full permutation group.\\

\n{\bf Definition 2.1}. [2;sec.3]. {\em For monomials

\[u=x^{j(1)}_{i(1)}x^{j(2)}_{i(\ell)}\cdots x^{j(\ell)}_{i(\ell)}\;\;\;\;\;
v=x^{j^{\prime}(1)}_{i^{\prime}(1)}x^{j^{\prime}(2)}_{i^{\prime}(2)}\cdots
x^{j^{\prime}(m)}_{i^{\prime}(m)}\]

\n in the free ring $P$, set $u>v$ when
\begin{list}%
{(\roman{rom})}{\usecounter{rom}\setlength{\rightmargin}{\leftmargin}}
\item deg $u>$ deg $v$, or
\item deg $u=$ deg $v$, the exponent sequences of $u$ and $v$ are the
  same (so $\ell=m$), and ``the first nonzero difference in subscripts
  $i-i^{\prime}$ is less than zero'', or
\item deg $u$=deg $v$, the exponent sequences of $u$ and $v$ are
  different, and ``the first nonzero difference in exponents
  $j-j^{\prime}$ is greater than zero''.  (Quoted phrases from Wolf
  [2])
\end{list}}

\n So monomials of the same degree are ordered lexicographically by
their complexions when they have the same exponent sequence, and
anti-lexicographically by their exponent sequences when these are
different.  Hence, for example, when $n=4$,
$x_{1}x_{2}^{2}x_{3}>x_{1}x^{2}_{2}x_{4}$ and
$x_{1}x^{2}_{4}>x_{1}x_{2}x_{3}$.  One checks directly that
Definition 2.1 provides a total ordering of the set of monomials.\\

\n Following Bergman and Cohn [1] we denote by $Q_{0}$ the set of
monomials that are maximal in their $G$-orbits.  Given a monomial $u$,
the $n$ (different) monomials in its $G$-orbit (equivalently, the $n$
monomials in its orbit polynomial $\overline{u}$) all have the same
exponent sequence, hence are ordered lexicographically by their
complexions according to 2.1{\em(ii)}.  Thus $Q_{0}$ {\em consists exactly
  of those monomials that begin with a power of $x_{1}$}.  Given
monomials $u,v$ in $Q_{0}$, define $u\cdot v$ to be the largest
monomial amoung the $n^{2}$ monomials in the (ordinary) product
$\overline{u}$ $\overline{v}$ of their orbit polynomials.  Define also
$1\cdot u=u\cdot 1=u$.  The product $\overline{u}$ $\overline{v}$ is
invariant, thus a sum (all coefficients 1) of orbit polynomials. Then
$u\cdot v$ is the largest amoung the maximal monomials of these orbit
polynomials.  This shows that $u\cdot v\in Q_{0}$, and, upon further
analysis; provides a simple procedure for computing $u\cdot v$: {\em
  the monomial $u\cdot v$ is the ordinary product of the largest
  monomial $y$ in $\overline{u}$ multiplied on the right by the (only)
  monomial in $\overline{v}$ whose first indeterminate equals the last
  one in $y$.} For example, $\overline{x}_{1}=x_{1}+x_{2}+\cdots
+x_{n}$ and $\overline{x_{1}x_{2}}=x_{1}x_{2}+x_{2}x_{3}+\cdots
+x_{n}x_{1}$, so $(x_{1})\cdot(x_{1}x_{2})=x^{2}_{1}x_{2}$.\\

\n {\bf Theorem 2.2,} {\em Under the product ``$\cdot$'',
  $Q_{0}\cup\{1\}$ is an associative semigroup whose atoms are
  precisely those monomials in $Q_{0}$ whose indeterminates occur only
  to the first power.  $Q_{0}\cup\{1\}$ is freely generated by its
  atoms.}\\

An ``atom'' in the semigroup $Q_{0}\cup\{1\}$ is a monomial that
cannot be written as the product of two monomials.  The product
$u\cdot v$ of two monomials $u,v$ in $Q_{0}$ computed as described
directly before the theorem, will always contain the square (at least)
of an indeterminate.  Hence a monomial whose indeterminates occur only
to the first power cannot be written as a product, so is an atom.
Conversely, suppose that $u\in Q_{0}$ is an atom.  If $u$ contained an
indeterminate to a power greater than one, then $u$ can be ``broken''
at that point, and written as product.  For example, when
$n=3,\;\;x_{1}x^{4}_{3}x^{2}_{1}x_{2}=(x_{1}x^{3}_{3})\cdot(x_{1}x^{2}_{2}x_{3})$.
So, if, $u$ is an atom, then its indeterminates occur only to the
first power.\\

\n ``Freely generated by its atoms'' means that every monomial in
$Q_{0}$ is a product of atoms, and that representation is unique: if
$a_{i},b_{j}$ are atoms, and $a_{1}\cdot a_{2}\cdot ... \cdot
a_{k}=b_{1}\cdot b_{2}\cdot ... \cdot b_{\ell}$, then $k=\ell$, and
$a_{i}=b_{i},i=1,2,...,k$. Theorem 2.2 corresponds to Proposition 2.1
in [1].\\

A listing of atoms would begin\\

\n deg 1 {\makebox{\hspace*{.4in}}} $x_{1}$\\
\n deg 2 {\makebox{\hspace*{.4in}}}
$x_{1}x_{2},x_{1}x_{3},...,x_{1}x_{n}$\\
\n deg 3 {\makebox{\hspace*{.4in}}}
$x_{1}x_{2}x_{1},x_{1}x_{2}x_{3},...,x_{1}x_{2}x_{n},...,x_{1}x_{n}x_{1},...,x_{1}x_{n}x_{n-1}$\\

\n There are $(n-1)^{d-1}$ atoms of degree $d$. Here are some examples
of factorizations into atoms.  Write $(x_{1})^{a}$ for $(x_{1})\cdot
(x_{1})\cdot ... \cdot (x_{1})$ (a times).

\[\begin{array}{l}x^{a}_{1}x^{b}_{i}=(x_{1})^{a-1}\cdot (x_{1}x_{i})\cdot
(x_{1})^{b-1}\;\;\;i\neq 1\\

x_{1}^{a}x^{b}_{i}x^{c}_{j}=(x_{1})^{a-1}\cdot (x_{1}x_{i})\cdot
(x_{1})^{b-2}\cdot (x_{1}x_{k})\cdot (x_{1})^{c-1}\;\;\;\;;i\neq
j\;\;i\neq 1
\end{array}\]

\n Here $b\geq 2$ and $k=j^{h}$ where $h=g^{-1}$ and $g\in G$
determined by $1^{g}=i$.

\[x^{a}_{1}x_{i}x^{c}_{j}=(x_{1})^{a-1}\cdot(x_{1}x_{i}x_{j})\cdot(x_{1})^{c-1}\;\;\;\;i\neq
j\;\;i\neq 1\]
\newpage

\n{\bf Theorem 2.3.} {\em Every polynomial in the free ring
  $P=K<x_{1},x_{2},...,x_{n}>$ that is invariant under the action of
the group $G=<(1 2 ...n)>$ of circular permutations is uniquely a
polynomial in orbit polynomials $\overline a_{i}$ where the $a_{i}$
are atoms in $Q_{0}$.}\\

This result corresponds to Theorem 1 and Corollary 1 in [1].  Theorem
2.3 is implemented exactly as in the commutative case.  Given a
$G$-invariant polynomial, it is uniquely a linear combination of orbit
polynomials, so it is enough to deal with orbit polynomials.  Let
$\overline{u}$ be an orbit polynomail where the monomial $u$ is chosen
to be maximal in its orbit, so $u\in Q_{0}$.  Factor $u$ into a
product of atoms:

\[u=a_{1}\cdot a_{2}\cdot ...\cdot a_{k}\; ,\; a_{i}\; \text{atoms in}\;
Q_{0},\; \text{repetitions possible}\]

\n As an ordinary product of monomials, we have

\[u=a_{1}a_{2}^{g}...a^{h}_{k}\;\;\text{where}\;g.....h\in G.\]

\n Now the product $\overline{a}_{1}\overline{a}_{2}\cdots
\overline{a}_{k}$ of the orbit polynomials $\overline{a}_{i}$ is
itself invariant, so is a sum of orbit polynomials.  Amoung the
maximal monomials of these, $u$ is biggest. Hence
$\overline{u}-\overline{a}_{1}...\overline{a}_{k}$ is a sum of orbit
polynomials whose maximal monomials are all smaller than $u$.  Repeat
the procedure with the largest amoung these.  The process will end in
a finite number of steps, leaving $\overline{u}$ expressed as a
polynomial in the $\overline{a}_{i}$.  Here are a couple of examples
of representations of orbit polynomials as polynomials in the
$\overline{a_{i}}$.  I have taken $n=3$ for simplicity.

\[\begin{array}{lcl}
\overline{x^{2}_{1}}=x^{2}_{1}+x^{2}_{2}+x^{2}_{3}&=&(\overline{x_{1}})^{2}-(\overline{x_{1}x_{2}}+\overline{x_{1}x_{3}})\\
\\
\overline{x^{3}_{1}}=x^{3}_{1}+x^{3}_{2}+x^{3}_{3} &=
&(\overline{x_{1}})^{3}+(\overline{x_{1}x_{2}x_{1}}+\overline{x_{1}x_{2}x_{3}}+\overline{x_{1}x_{3}x_{1}}+\overline{x_{1}x_{3}x_{2}})\\
\\
 &- &
 ((\overline{x_{1}})(\overline{x_{1}x_{2}})+(\overline{x_{1}x_{2}})(\overline{x_{1})}+(\overline{x_{1}})(\overline{x_{1}}x_{3})+(\overline{x_{1}x_{3}})(\overline{x_{1}}))
\end{array}\]

\n So we have an explicit description of the polynomials in the free
ring $P$ that are fixed under $G$: each is uniquely a polynomial in
the orbit polynomials $\overline{a}$ where the $a$'s are monomials that
begin with $x_{1}$ and whose indeterminates occur only to the first
power.  Each such polynomial, considered now in $R^{\prime}=P/I$,
  remains invariant there (because $I^{g}=I$ for all $g$ in $G$), and
  has the same representation as a polynomial in the $\overline{a}$,
  but no longer unique.  Moreover, all invariant polynomials in $R$
  arise this way.\\

\n{\bf Theorem 2.4} {\em Every $K$-coefficient polynomial in $R$ fixed
  by $G$ is a $K$-coefficient polynomial in the orbit polynomials
  $\overline{a}=\sum(a^{g}:\;g\in G)$, where the monomials $a$ begin
  with $x_{1}$ and have all their indeterminates to the first
power.}\\

\n {\bf Proof.} Given $p\in R$ fixed by $G$ consider the same
polynomial in the free ring $P$.  It may or may not be invariant
there, but $\overline{p}=\sum(p^{g}:g\in G)$ is.  Apply $\varphi:
P\longrightarrow
P/I:\;\;\varphi(\overline{p})=\sum\varphi(p^{g})=\sum\varphi(p)^{g}=np,$
so $p=(1/n)\varphi(\overline{p})$ expresses $p$ in the form
desired.$\diamond$\\

\n For example, when $n=3$, $p=x_{1}x_{2}-x_{2}x_{1}$ is invariant in $R$
but not in $P$.  Write
$\overline{p}=\overline{x_{1}x_{2}}-\overline{x_{1}x_{3}}$, whence
$x_{1}x_{2}-x_{2}x_{1}=(1/3)(\overline{x_{1}x_{2}}-\overline{x_{1}x_{3}})$
in $R$.\\

\n The $\sigma_{k}$ are invariant in $R$, so they fall under Theorem
2.4.  In this case the $a's$ that appear have increasing complexions,
and the representation is reasonably explicit:\\

\n {\bf Corollary 2.5.} {\em The $\sigma_{k}(x_{1},x_{2},...,x_{n}),
k=1,2,...,n,$ can be written in $R$ as a linear combination of orbit
polynomials with positive integer coefficients as follows

\[n\sigma_{k}=\sum(\alpha_{j}\overline{x_{1}x_{i(2)}\cdots
x_{i(k)}}:\;i(\cdot)\in\{2,3,...,n\},\;1<i(2)<\cdots<i(k))\]

\n where the positive integers $\alpha_{j}$ are easily computed.}\\

\n The result follows from formula (2) of Sec 1; I omit the details.
Here is a sample calculation for $n=4$.

\[\begin{array}{rcl}
\sigma_{1} & = & \overline{x_{1}}\\
4\sigma_{2} & = &
3\overline{x_{1}x_{2}}+ 2\overline{x_{1}x_{3}}+\overline{x_{1}x_{4}}\\
4\sigma_{3} &=&
2\overline{x_{1}x_{2}x_{3}}+\overline{x_{1}x_{2}x_{4}}+\overline{x_{1}x_{3}x_{4}}\\
4\sigma_{4} &=& \overline{x_{1}x_{2}x_{3}x_{4}}
\end{array}\]

\n So, just as in the free ring $P$, every invariant polynomial in
$R$ is a polynomial in the orbit polynomials $\overline{a}$, where the
monomials $a$ begin with $x_{1}$ and have all their indeterminates to
the first power.  In $P$, these orbit polynomials are algebraically
independent, and there are infinitely many of them.  In $R$, the
central issue regarding its fixed ring is to select from these
``residual atoms'' a smaller generating set, hopefully finite.  When
$n=3$ it turns out that one will do (next section).\\

\vspace*{.2in}
\n{\bf THE CASE $n=3$.}\\

\n We assume $n=3$ throughout this section.  So here we are dealing with
$R^{\prime}=P/I$ where $I$ is the ideal of the free ring
$P=K<x_{1},x_{2},x_{3}>$ generated by the four polynomials

\[\begin{array}{lcl}
A & = & \sigma_{2}-\sigma_{2}^{g(1)}=[1,2]+[1,3]\\
B & = & \sigma_{2}^{g(1)}-\sigma_{2}^{g(2)}=[2,3]+[2,1]\\
C & = &
\sigma_{3}-\sigma_{3}^{g(1)}=x_{1}x_{2}x_{3}-x_{2}x_{3}x_{1}=[1,23]=x_{2}[1,3]+[1,2]x_{3}\\
D & = &
\sigma_{3}^{g(1)}-\sigma_{3}^{g(2)}=x_{2}x_{3}x_{1}-x_{3}x_{1}x_{2}=[2,31]=x_{3}[2,1]+[2,3]x_{1}
\end{array}\]

\n Our augmented ring $R$ contains the field
$W=K(\sigma_{1},\sigma_{2},\sigma_{3})$ within its center, the
elementary polynomials
$\sigma_{1}=x_{1}+x_{2}+x_{3},\;\;\sigma_{2}=x_{1}x_{2}+x_{1}x_{3}+x_{2}x_{3},\;\;\sigma_{3}=x_{1}x_{2}x_{3}$
being algebraically independent over $K$ and invariant in $R$ under
the group $G=<(123)>$.  \\

\n{\bf Definition 3.1.} {\em Denote by $c$ the common value
  $[1,2]=[2,3]=[3,1]$ which is nonzero and invariant (refer to Theorem
  1.3).}\\
\newpage

\n{\bf Lemma 3.2} {\em In the ring $R$, 
\begin{list}%
{(\roman{rom})}{\usecounter{rom}\setlength{\rightmargin}{\leftmargin}}
\item $c^{3},$ and the two orbit polynomials
  $\overline{x_{1}x_{2}x_{1}}$ and $\overline{x_{1}x_{3}x_{1}}$, are
each nonzero, central, and invariant.  None belong to $W$.
\item Neither $c$ nor the orbit polynomial $\overline{x_{1}x_{2}}$
  commutes with $x_{1},x_{2},$ or $x_{3}$. Each has degree 3 over $W[c^{3}]$.
\end{list}}

\n{\bf Proof.} (i) The relations $A=B=C=D=0$ hold in $R$.  $C=0$ gives
$x_{2}c=cx_{3},$ and $D=0\;\;x_{3}c=cx_{1}$.  Applying $G$ we get
$x_{1}c=cx_{2}$.  These three in turn yield $x_{1}c^{2}=c^{2}x_{3},\;\;
x_{2}c^{2}=c^{2}x_{1},$ and $x_{3}c^{2}=c^{2}x_{2}$. From these we get
$x_{i}c^{3}=c^{3}x_{i},\;\;i=1,2,3,$ thus $c^{3}$ is central.  As for
$c^{3}\neq 0$, the four relations above yield all identities that hold
in $R$.  We have just shown that $C=0$ and $D=0$ are equivalent to
$x_{i}c=cx_{i+1}(4\equiv 1)$.  The relations $A=0,\;B=0$ are
equivalent to $[1,2]=[2,3]=[3,1]$.  Clearly $c^{3}=0$ is not a
consequence of these. Obviously $c^{3}$ is invariant because $c$ is.
As for $c^{3}\not\in W$, suppose that it were.  Then $c^{3}$ would
equal a rational function in $\sigma_{1},\sigma_{2},\sigma_{3}$ whose
numerator would not be zero because $c^{3}\neq 0$. That same equality
will hold when the variables $x_{i}$ are allowed to commute (Theorem
1.7(i)). Then $c^{3}=0$ but the right side
remains unchanged which contradicts the algebraic independence of the
$\sigma_{i}$. Hence $c^{3}\not\in W.$ As for the orbit polynomials
$\overline{x_{1}x_{2}x_{1}}$ and $\overline{x_{1}x_{3}x_{1}}$ they are
obviously invariant and neither is zero because they cannot be brought
to zero by allowing the $x_{i}$ to commute.  And neither can lie in
$W$ because, for example, were $\overline{x_{1}x_{2}x_{1}},\in W$ then
  it would equal a rational function in the $\sigma_{k}$.  But,
  allowing the variables to commute, we would have
  $x^{2}_{1}x_{2}+x^{2}_{2}x_{3}+x^{2}_{3}x_{1}$ equal to a quantity
  invariant under $S_{3}$, a contradiction. Finally, as for their centrality, we need only consider one by
virtue of the relation
$\overline{x_{1}x_{2}x_{1}}+\overline{x_{1}x_{3}x_{1}}=\sigma_{1}\sigma_{2}-3\sigma_{3}$.
    To check this, argue as follows:
    $x_{1}x_{2}x_{1}+x_{1}x_{3}x_{1}=x_{1}(x_{2}+x_{3})x_{1}=x_{1}(\sigma_{1}-x_{1})x_{1}=\sigma_{1}x_{1}^{2}-x^{3}_{1}=\sigma_{2}x_{1}-\sigma_{3}$, the last step because $x_{1}$ is a root of $y^{3}-\sigma_{1}y^{2}+\sigma_{2}y-\sigma_{3}$. Then pass to the orbit polynomials.  We shall work with $\overline{x_{1}x_{2}x_{1}}$.  Because it is invariant, to prove that it is central, it is enough to prove that it commutes with $x_{1}$; i.e. that $[x_{1},\overline{x_{1}x_{2}x_{1}}]=0$. In the following computation remember that, in $R$, $\sigma_{3}$ is central and invariant so $\sigma_{3}=x_{1}x_{2}x_{3}=x_{2}x_{3}x_{1}=x_{3}x_{1}x_{2}$. We shall use the simplifying notation $[x_{1},\overline{x_{1}x_{2}x_{1}}]=[1,\overline{121}].$

\[\begin{array}{lcl}
[1,\overline{121}] & = & [1,(12)1]+[1,(23)2]+[1,(31)3]\\
                   & = & [1,12]1+23[1,2]+[1,23]2+31[1,3]+[1,31]3\\
                   & = & 1[1,2]1+23[1,2]+31[1,3]+[1,3]13 \ \ \ ([1,23]=0)\\
                   & = &
                   x_{1}cx_{1}+x_{2}x_{3}c-x_{3}x_{1}c-cx_{1}x_{3}\\
                   & = & x_{1}x_{3}
                   c+x_{2}x_{3}c-x_{3}x_{1}c-x_{3}x_{2}c\\
                   & = & ([1,3]+[2,3])c=(-c+c) c=0 c=0.
\end{array}\]

\n Hence $\overline{x_{1}x_{2}x_{1}}$ and $\overline{x_{1}x_{3}x_{1}}$
are central.  The fact that $c^{3},$ which has degree 6, and the
cubics $\overline{x_{1}x_{2}x_{1}}$ and $\overline{x_{1}x_{3}x_{1}}$ are central contrasts with the
fact that every central quadratic polynomial in $R$ is a
$K$-coefficient polynomial in $\sigma_{1}$ and $\sigma_{2}$ (proof
omitted).

\n (ii) Because $\overline{x_{1}x_{2}}$ and $c$ are related by the
easily verified formula $\overline{x_{1}x_{2}}=c+\sigma_{2}$, it is
enough to deal with $c$.  And, as $c$ is invariant, we need only show
that it does not commute with $x_{1}$.  Note that if we knew that $R$
was a domain (Conjecture 1.9), the proof would be trivial: from
$x_{1}c=cx_{1}$ and $x_{1}c=cx_{2}$ we would have
$c(x_{1}-x_{2})=0$ a contradiction as neither is zero.  To avoid
relying on the truth of Conjecture 1.9 at this point, I supply a
considerably longer proof which goes as follows. Use $c=[2,3]$ and
compute $[1,c]=[1,[2,3]]=-x_{1}x_{3}x_{2}+x_{3}x_{2}x_{1}$. So it is a
matter of showing that $x_{1}x_{3}x_{2}-x_{3}x_{2}x_{1}\neq 0$ in $R$.
I shall prove this by showing that, in the free ring
$P,\;\;x_{1}x_{3}x_{2}-x_{3}x_{2}x_{1}\not\in I$. A homogeneous cubic
will belong to $I$ if and only if it can be written as a $K$-linear
combination of the following four polynomials: $qAr,\;\;sBt,\;\;uCv.$
amd $yDz,$ where the eight monomials involved must satisfy $\deg
q+\deg r=\deg s+\deg t=1$ and $\deg u=\deg v=\deg y=\deg z=0$.  So we
have the following six expressions:

\[\begin{array}{ll}
(\alpha_{1}1+\alpha_{2}2+\alpha_{3}3) & (12-21+13-31)\\
                                    &(12-21+13-31)(\beta_{1}1+\beta_{2}2+\beta_{3}3)\\
(\gamma_{1}1+\gamma_{2}2+\gamma_{3}3 & (23-32+21-12)\\
                                     &
                                     (23-32+21-12)(\delta_{1}1+\delta_{2}2+\delta_{3}3)\\
                                     &\mu(123-231)\\
                                     &\sigma(231-312)
\end{array}\]

\n where I have written 12 for $x_{1}x_{2},$ etc, and the greek
letters stand for elements of $K$. There are 27 monomials of degree 3.
We may disregard the three that have exponent sequence 3.  I have
listed below the remaining 24 in decreasing order.

\[\begin{array}{c|c}
Exponent & \\
sequence & Monomials\\ \hline
21       & 1^{2}2>1^{2}3>2^{2}1>2^{2}3>3^{2}1>3^{2}2\\
\vee        & \\
12       & 12^{2}>13^{2}>21^{2}>23^{2}>31^{2}>32^{2}\\
\vee        & \\
111      & 121>123>131>132>212>213>231>232\\
         & >312>313>321>323
\end{array}\]

\n Only 132-321 appears on the right side of our equation, so all
other monomials must have coefficient zero.  Beginning with the
monomials with exponent sequence 21 we get

\[\begin{array}{c}
(\alpha_{1}-\gamma_{1})1^{2}2+\alpha_{1}1^{2}3+(-\alpha_{2}+\gamma_{2})2^{2}1-\gamma_{2}2^{2}3\\
+(-\alpha_{3})3^{2}1+(-\gamma_{3})3^{2}2
\end{array}\]

\n which implies
$\alpha_{1}=\alpha_{2}=\alpha_{3}=\gamma_{1}=\gamma_{2}=\gamma_{3}=0$.
Working similarly with exponent sequence 12 we get
$\beta_{1}=\beta_{2}=\beta_{3}=\delta_{1}=\delta_{2}=\delta_{3}=0$. So
we are left with

\[\mu(123-231)+\sigma(231-312)=132-321\]

\n which is clearly impossible.  Thus our result: $c$, so also
$\overline{x_{1}x_{2}}=c+\sigma_{2},$ commutes with neither
$x_{1},x_{2},$ nor $x_{3}$.  As for the $Z$-coefficient cubic
satisfied by $\overline{x_{1}x_{2}}$, just cube the equation $\overline{x_{1}x_{2}}-\sigma_{2}=c$ to get
$(\overline{x_{1}x_{2}})^{3}-3\sigma_{2}(\overline{x_{1}x_{2}})^{2}+3\sigma_{2}(\overline{x_{1}x_{2}})-(\sigma_{2}^{3}+c^{3})=0$.
$\diamond$\\

\n{\bf Definition 3.3.} {\em Let
  $Z=W[c^{3},\overline{x_{1}x_{2}x_{1}}]$ which is a central subring of
    $R$ pointwise fixed by $G$. (The bracket [ ] here denotes
    adjunction, not the commutator.) Denote by $S$ the subring of $R$
    consisting of all polynomials fixed by $G$.}\\

\newpage
\n{\bf Theorem 3.4.} {\em Every $p\in S$ has the form
  $z_{0}+z_{1}c+z_{2}c^{2}$ for some $z_{i}\in Z$.  Accordingly
\begin{list}%
{(\roman{rom})}{\usecounter{rom}\setlength{\rightmargin}{\leftmargin}}
\item $S$ is a commutative degree 3 extension of $Z$, and
\item For every $p\in R,\;\;\overline{p}=\sum(p^{g}:g\in G)$ is either
  central or is a root of a cubic, coefficients in $Z$.
\end{list}}

\n{\bf Proof.} In this proof we shall refer to monomials in $R$ that
begin with $x_{1}$ and have all their indeterminates to the first
power as ``residual atoms''.  These are not atoms in any technical
sense - this is just a handy notation.  In Theorem 2.4 we have shown
that every $K$-coefficient polynomial in $R$ that is invariant under
$G$ is a $K$-coefficient polynomial in the orbit polynomials of
residual atoms.  So it is enough to prove that every orbit polynomial
of a residual atom is a polynomial in $\overline{x_{1}x_{2}},$
coefficients in $Z$.  We begin by establishing this for the residual
atoms of degree $\leq 3$:

\[\begin{array}{l}
\overline{x_{1}} =  \sigma_{1}\in Z\\
\overline{x_{1}x_{3}}  = 
3\sigma_{2}-2\overline{x_{1}x_{2}}=3\sigma_{2}-2(c+\sigma_{2})=\sigma_{2}-2c\not\in
Z\\
\overline{x_{1}x_{2}x_{1}}\in Z \\
\overline{x_{1}x_{2}x_{3}}  =  3\sigma_{3}\in Z\\
\overline{x_{1}x_{3}x_{1}}  = 
\sigma_{1}\sigma_{2}-3\sigma_{3}-\overline{x_{1}x_{2}x_{1}}\in Z\\
\overline{x_{1}x_{3}x_{2}} = 
  \sigma_{1}\sigma_{2}+3\sigma_{3}-(\overline{x_{1}x_{2}})\sigma_{1}=3\sigma_{3}-c\sigma_{1}\not\in Z\\
\end{array}\]

\n Relations 1,3, and 4 are evident; number 5 was done in the proof of
Lemma 3.2.  We begin with 2:

\[\begin{array}{lcl}
3\sigma_{2}-2\overline{x_{1}x_{2}} & = &
\sigma_{2}+2(\sigma_{2}-\overline{x_{1}x_{2}})=\sigma_{2}-2c=\sigma_{2}-c-c\\
 & = & \sigma_{2}-(x_{1}x_{2}-x_{2}x_{1})-(x_{2}x_{3}-x_{3}x_{2})\\
 & = &
 x_{1}x_{2}+x_{1}x_{3}+x_{2}x_{3}-x_{1}x_{2}+x_{2}x_{1}-x_{2}x_{3}+x_{3}x_{2}\\
 & = & x_{1}x_{3}+x_{2}x_{1}+x_{3}x_{2}=\overline{x_{1}x_{3}}.
\end{array}\]

\n Finally number 6, which uses
$c=x_{2}x_{3}-x_{3}x_{2}=\overline{x_{1}x_{2}}-\sigma_{2}$:

\[\begin{array}{lcl}
x_{1}x_{3}x_{2} & = &
x_{1}(x_{2}x_{3}-c)=x_{1}(x_{2}x_{3}+\sigma_{2}-\overline{x_{1}x_{2}})\\
& = &
x_{1}x_{2}x_{3}+x_{1}\sigma_{2}-x_{1}(\overline{x_{1}x_{2}}),\;\;whence\\
\overline{x_{1}x_{3}x_{2}} & = &
\overline{x_{1}x_{2}x_{3}}+\overline{x_{1}}\sigma_{2}-\overline{x_{1}}(\overline{x_{1}x_{2}})\\
& = &
3\sigma_{3}+\sigma_{1}\sigma_{2}-\sigma_{1}(\overline{x_{1}x_{2}}).
\end{array}\]

\n Let us call a residual atom {\em a reducible} if $\overline{a}$ is
a $Z$-coefficient polynomial in orbit polynomials of residual atoms of
degree strictly less than that of $a$. We shall prove that {\em every
  residual atom of degree $\geq 4$ is reducible} which will establish
Theorem 3.4.  A residual atom $a$ of degree $\geq 4$ will begin
121..., 123..., 131..., or 132... . If $a=121u$ then the monomial $u$
must begin with either 2 or 3.  In the first case use
$12=\sigma_{2}-13-23$ to get
$a=1212...=(\sigma_{2}-13-23)12...=(12)\sigma_{2}-1312...-2312...$.
In $R,\sigma_{3}=312$, so we get
$\overline{a}=(\overline{12})\sigma_{2}-\sigma_{3}(\overline{1...})-\sigma_{3}(\overline{2...})$.
In the second case
$a=1213...=1(-c+12)3...=-1c3...-1(123)...=-c23...-\sigma_{3}(1...)=-(\overline{x_{1}x_{2}}-\sigma_{2})23...-\sigma_{3}(1...)$
  So
  $\overline{a}=-(\overline{x_{1}x_{2}})(\overline{23...})-\sigma_{2}(\overline{23...})-\sigma_{3}(\overline{1...})$.  If $a=123u$ then $\overline{a}=\sigma_{3}\overline{u}$.  If $a=131u$ and $u$ begins with 2, then $a=1312...=\sigma_{3}(1...)$. so $\overline{a}=\sigma_{3}(\overline{1...})$. If $u$ begins with 3 so $a=1313...$, then use $\sigma_{2}=23+21+31$ to write $a=1(\sigma_{2}-23-21)3...=\sigma_{2}(13...)-(123)(3...)-1213....$ Then $\overline{a}=\sigma_{2}(\overline{13...})-\sigma_{3}(\overline{3...})-\overline{1213}...$ and the last term begins with $121$ so is reducible by our first case.  Finally, if $a=132u$ then we write $13=31-c=31-\overline{12}+\sigma_{2}$ to get $a=(31-\overline{12}+\sigma_{2})2u=312u-(\overline{12})(2u)+\sigma_{2}(2u)$ whence $\overline{a}=\sigma_{3}\overline{u}-(\overline{12})(\overline{2u})+\sigma_{2}(\overline{2u})$. $\diamond$\\

\n{\bf Theorem 3.5.} {\em If $R$ is a domain, then
\begin{list}%
{(\roman{rom})}{\usecounter{rom}\setlength{\rightmargin}{\leftmargin}}
\item its center is $Z$,
\item both $c$ and $\overline{x_{1}x_{2}x_{1}}$ are algebraically
  independent transcendentals over $W$, and
\item when $R$ is enlarged as described below, every nonzero invariant
  polynomial is invertible in $R$.
\end{list}}

\n {\bf Proof.} (i) If $p$ is central it commutes with $c$ so $cp=pc$.
Multiplication on the right by $c$ has the effect of applying
$g=(123),\;pc=cp^{g}$, whence $c(p-p^{g})=0$ so that $p$ is invariant.
Then $p=z_{0}+z_{1}c+z_{2}c^{2}$ must commute with $x_{1}$ which
forces $z_{1}=z_{2}=0$\\
(ii) We shall first prove that $c$ and $\overline{x_{1}x_{2}x_{1}}$
  are individually transcendental over $W$.  Suppose that $c$ were a
  root of a monic polynomial, $\sum a_{j}c^{j}=0$, where the $a_{j}\in
  W$.  If this polynomial has no constant term, we may factor out a
  power of $c$ to get $c^{k}\sum b_{j}c^{j}=0$ where the latter
  polynomial now has a nonzero constant term. It must itself be zero
  as we are working in a domain and $c\neq 0$.  Now allow the $x_{i}$
  to commute.  All $c$'s will vanish leaving our constant term equal
  to zero contradicting the algebraic independence of the ``elementary
  symmetric functions'' in commuting variables. The proof for
  $\overline{x_{1}x_{2}x_{1}}$ is similiar except that, when the
  $x_{i}$ are allowed to commute, we reach the contradiction that a
  polynomial in $x^{2}_{1}x_{2}+x^{2}_{2}x_{3}+x^{2}_{3}x_{1}$, which
is not invariant under $S_{3}$, equals a rational function in the
$\sigma_{k}$, which is.  Finally, the algebraic independence of $c$
and $\overline{x_{1}x_{2}x_{1}}$.  Set $d=\overline{x_{1}x_{2}x_{1}}$
for convenience.  Let $\sum a_{i_{j}}c^{i}d^{j}$ be a monic polynomial
in which $c$ and $d$ both actually appear, $a_{ij}\in W$.  If there is
no constant term we may factor out $c^{r}d^{s}$ (one of the integers
$r$ or $s$ may be zero) to get $c^{r}d^{s}\sum b_{ij} c^{i}d^{j}$ where
the latter polynomial does have a nonzero constant term.  If the
original polynomial is zero, so is the latter.  Then argue as before.\\
(iii)Owing to the fact that $c^{3}$ and $\overline{x_{1}x_{2}x_{1}}$
  are central algebraically independent transcendentals over $W$, we
  may enlarge $R$ using the same method described in the paragraph
  following Definition 1.8 to a ring that has the field of fractions
  $W(c^{3},\overline{x_{1}x_{2}x_{1}})$ as its center.  Let us use the
  same symbol $R$ for this enlarged ring, and the same symbol $Z$ for
  its center $W(c^{3},\overline{x_{1}x_{2}x_{1}})$ which we note has
  transcendence degree 5 over $K$.  Then, in this $R$, every nonzero
  invariant polynomial either lies in the center $Z$, in which case it
  is invertible, or satisfies a cubic with nonzero constant term that
  lies in $Z$, thus invertible in this case as well. $\diamond$\\

\n I have not been able to answer the natural follow-up question to
Theorem 3.5: is $R$ itself a field?  More to the point, does this
method of construction (forming the quotient of the free ring
$K<x_{1},x_{2},...,x_{n}>$ by the ideal $I$ generated by the
$[x_{i},\sigma_{k}]$) provide a way to construct finite-dimensional
fields? If so, do said fields have any novel properties?  All this is
related to the long standing open question: is every field of degree 5
cyclic?

\end{document}